\documentclass[11pt]{article} 
\usepackage{graphicx}
\usepackage{latexsym}
\usepackage{amsmath}
\usepackage{amssymb}
\usepackage{amsthm}
\usepackage{verbatim}
\usepackage{tikz}

\newcommand{\la}{\langle}
\newcommand{\ra}{\rangle}

\title{Single-delete Nim}
\author{Masato Shinoda\thanks{Division of Natural Sciences, Nara Women's University. {\tt shinoda@cc.nara-wu.ac.jp}}}

\date{\today}

\begin{document}

\theoremstyle{definition} 
\newtheorem{theorem}{Theorem}[section]
\newtheorem{definition}[theorem]{Definition}
\newtheorem{lemma}[theorem]{Lemma}
\newtheorem{proposition}[theorem]{Proposition}
\newtheorem{corollary}[theorem]{Corollary}
\newtheorem{example}[theorem]{Example}
\newtheorem{remark}[theorem]{Remark}
\newtheorem{conjecture}[theorem]{Conjecture}

\maketitle              
\begin{abstract}
The classic game of Nim has been well-known for many years, inspiring numerous variations. One such variant is Delete Nim, where players take turns eliminating one pile of stones and splitting the remaining pile into two smaller piles. 
In this paper we generalize the game to include the case of $n$ piles. On each turn, a player eliminates one pile and splits one of the remaining piles into two smaller piles. We specifically analyze the case where $n=4$, deriving the conditions for a winning strategy.

\end{abstract}
\section{Introduction and Main Result}
\subsection{Introduction}

Nim is a well-known game that has been studied extensively over the years. In the most popular version of Nim, there are three piles of stones, and two players take turns. On their turn, a player selects one pile and removes any positive number of stones from that pile. The players continue taking turns, and the player who removes the last stone wins the game. This game can also be played with more than three piles following the same rules.

Nim is a two-player, zero-sum, finite, deterministic, and perfect information game with no possibility of a draw. Therefore, every configuration can be classified as either an N-position, where the player whose turn it is can force a win, or a P-position, where the player whose turn it is will inevitably lose if the opponent plays optimally. 
The classification of N-positions and P-positions in the aforementioned version of Nim was demonstrated by Bouton \cite{bou02} in 1902. 
It became known that the conditions for the winning strategy in this game include mathematically intriguing structures.

Several variants of Nim have been proposed by modifying the rules of legal moves, such as Moore's game \cite{moo10} and Wythoff's game \cite{wyt07}. The winning strategies of these variants have been studied extensively, and their mathematical formulations have attracted significant interest. In this study, we generalize Delete Nim, a game introduced by Abuku and Suetsugu \cite{abu21}, where players eliminate one pile and then split one of the remaining piles into two smaller piles in one move. We then investigate the winning strategy for this generalized version.

\subsection{Delete Nim}

In this section, we first describe the rules of Delete Nim, which form the basis of this study. There are two variations of the Delete Nim rules provided in \cite{abu21} and \cite{sta08}. Although the rules differ in presentation, they are equivalent with respect to legal moves. In this paper, we adopt the latter rules\footnote{In \cite{abu21}, the rule described here is referred to as the ``Variant of Delete Nim'' (VDN).}.

In Delete Nim, there are two piles of stones. On each turn, players perform the following two actions sequentially (we refer to these combined actions as {\it a move}),  then pass the turn to their opponent.
\begin{itemize}
\item Select one of the piles and remove all the stones from it, eliminating the pile.
\item Split the remaining pile into two new piles, ensuring that each new pile contains at least one stone.
\end{itemize}
A player unable to make this move on their turn loses the game. We denote a position in Delete Nim by $\la y,z \ra$, where $y$ and $z$ represent the number of stones in the two piles. 
Let $\mathbb{N}$ be the set of positive integers. Then, the set of all possible game positions (without specifying the current player's turn) is denoted by $G_{2}=\{\, \la y,z \ra \, | \, y,z\in\mathbb{N}\}$. 
At the position $\la 1,1 \ra$, the player whose turn it loses the game.

Similar to Nim, Delete Nim is a two-player, zero-sum, perfect information, deterministic, finite game with no possibility of a draw. Therefore, each position in $G_2$ can be classified as either an N-position or a P-position. The condition for determining this classification is as follows.
\begin{theorem}[\cite{abu21,sta08}]
In Delete Nim, a position $\la y,z \ra$ is a P-position if both $y$ and $z$ are odd; otherwise, it is an N-position.
\end{theorem}
This theorem is included in Proposition 4.1 in this paper. It should be noted that in \cite{abu21}, not only is the above classification provided, but the Sprague-Grundy number of each position is also determined.
In an N-position, a player can ensure winning by forcing their opponent into P-positions with optimal moves. 
For example, when the current position is $\la 4,3 \ra$, the player on turn can proceed as $\la 4,3 \ra \to \la 3,1\ra$ (and then the game proceeds as $\to \la 2,1 \ra \to \la 1,1\ra$ and ends) and secure winning.
Therefore, identifying the P-positions is equivalent to understanding a winning strategy and knowing the conditions for winning.

\subsection{Single-delete Nim}

In this section, we define an extended rule that generalizes Delete Nim to any number of piles. Let $n$ be an integer greater than or equal to 2.
\begin{definition}[Single-delete Nim] 
There are $n$ piles of stones. On each turn, the two players perform the following two actions sequentially, then pass the turn to their opponent.
\begin{itemize}
\item Select one of the piles and remove all the stones from it, eliminating the pile.
\item Select one of the remaining $n-1$ piles and split it into two new piles, ensuring that each new pile contains at least one stone. 
\end{itemize}
A player unable to make this move on their turn loses the game.
\end{definition}
In Single-delete Nim, it is important to note that the number of piles remains $n$ at the end of each turn.
\begin{remark}
We introduce the term ``Single-delete Nim'' in this paper. When generalizing the number of piles in Delete Nim, various extended rules can be considered depending on how piles are eliminated and split. The name ``Single-delete Nim''  explicitly refers to the rule of  eliminating only one pile. For more details, see Abuku et al. \cite{abu24}.
\end{remark}
The following result for Single-delete Nim with three piles is known from Sakai \cite{sak21}. 
To present the winning and losing conditions, we define the following notation: For $z\in\mathbb{N}$, if $z$ is divisible by $2^{k}$ but not divisible by $2^{k+1}$, we denote it as $v_{2}(z)=k$. When $z$ is odd, we define $v_{2}(z)=0$.

\begin{theorem} A necessary and sufficient condition for the position $\la x,y,z \ra$ in Single-delete Nim with three piles to be a P-position is
\begin{quote}
\quad ($\bigstar$) \qquad\qquad $v_{2}(x)=v_{2}(y)=v_{2}(z)$.
\end{quote}
\end{theorem}
A proof of Theorem 1.4 is provided in \cite{sak21} (pp.59-63).
However, as it is also relevant for the case of four piles, we restate it in Section 2.1 of this paper.

\subsection{Main Result}
As the main theorem of this paper, we present the necessary and sufficient condition to be a P-position in Single-delete Nim with four piles as Theorem 1.5. 
To state this theorem, we introduce notation related to the binary representation of positive integers. Let $I_{k}(z)$ represent the $k$-th digit from the right in the binary representation of $z\in\mathbb{N}$. In other words, if the quotient of $z$ divided by $2^{k-1}$ is even, then $I_{k}(z)=0$; if the quotient is odd, then $I_{k}(z)=1$.
We note that $v_{2}(z)=\min\{\, k \, | \, I_{k}(z)=1\}-1$.

\begin{theorem}
In Single-delete Nim with four piles, for a position $\la w,x,y,z \ra$,
let $a=v_{2}(w), \, b=v_{2}(x), \, c=v_{2}(y), \, d=v_{2}(z)$. 
When $a\leq b\leq c\leq d$, the necessary and sufficient condition for $\la w,x,y,z \ra$ to be a P-position
is that one of the following conditions (1) to (5) is satisfied.

\begin{enumerate}
\renewcommand{\theenumi}{(\arabic{enumi})}
\setlength{\leftskip}{8mm}
\item $a=b=c=d$. 
\item $a<b=c=d$ and \\
(2A) \ $I_{b+1}(w)=0$. 
\item $a<b<c=d$ and all of (3A)-(3C) are satisfied.\\ 
(3A) \ $I_{c+1}(w)=I_{c+1}(x)=0$, \\ 
(3B) \ $I_{k}(w)+I_{k}(x)\geq 1$ for $b+2\leq k \leq c$, \\ 
(3C) \ $I_{b+1}(w)=1$.
\item $a<b<c<d$ and all of (4A)-(4E) are satisfied.\\ 
(4A) \ $I_{d+1}(w)=I_{d+1}(x)=I_{d+1}(y)=0$,\\
(4B) \ $I_{j}(w)+I_{j}(x)+I_{j}(y)\geq 2$ for $c+2\leq j \leq d$, \\
(4C) \ $I_{c+1}(w)=I_{c+1}(x)=1$,\\ 
(4D) \ $I_{k}(w)+I_{k}(x)\geq 1$ for $b+2\leq k \leq c$,\\
(4E) \ $I_{b+1}(w)=1$.
\item $a<b<c<d$ and all of (5A)-(5F) are satisfied. \\
(5A) \ $I_{i}(w)+I_{i}(x)+I_{i}(y)+I_{i}(z)\in\{0,3,4\}$ for $i\geq d+2$,\\
(5B) \ $I_{d+1}(w)=I_{d+1}(x)=I_{d+1}(y)=1$,\\
(5C) \ $I_{j}(w)+I_{j}(x)+I_{j}(y)\geq 2$ for $c+2\leq j \leq d$, \\
(5D) \ $I_{c+1}(w)=I_{c+1}(x)=1$,\\ 
(5E) \  $I_{k}(w)+I_{k}(x)\geq 1$ for $b+2\leq k \leq c$, \\
(5F) \ $I_{b+1}(w)=1$.
\end{enumerate}
\end{theorem}

The organization of this paper is as follows. In Section 2, we present a proof of Theorem 1.4, which provides the classification of N-positions and P-positions for Single-delete Nim with three piles, along with several useful propositions. In Section 3, we give a proof of Theorem 1.5, which addresses the case with four piles. In Section 4, we offer a brief analysis of the scenario when the number of piles is five or more.

\section{Single-delete Nim with three piles}
In this section, we address a proof of Theorem 1.4 for the case with three piles.
To support this proof, we prepare the following propositions.
These are nearly self-evident facts, so the proofs are omitted here, but we present them collectively as they will be frequently used in the next chapter as well.

\begin{proposition} 
For $x,y,z\in\mathbb{N}$ satisfying $x+y=z$, the following facts (i) and (ii) hold.
\begin{enumerate}
\renewcommand{\theenumi}{(\roman{enumi})}
\setlength{\leftskip}{8mm}
\item If $v_{2}(x)=v_{2}(y)$, then $v_{2}(z)>v_{2}(x)$, 
\item If $v_{2}(x)\neq v_{2}(y)$, then $v_{2}(z)=\min\{ v_{2}(x),v_{2}(y) \}$. 
\end{enumerate}
\end{proposition}

\begin{proposition} 
When $z\in\mathbb{N}$ and $v_{2}(z)\geq 1$, for any non-negative integer $k<v_{2}(z)$, there exist $x,y\in\mathbb{N}$ such that  $x+y=z$ and $v_{2}(x)=v_{2}(y)=k$.
\end{proposition}


\medskip

{\bf Proof of Theorem 1.4.} Let $P$ denote the set of all positions that satisfy condition ($\bigstar$). This set is a subset of  $G_{3}=\{\, \la x,y,z \ra \, | \, x,y,z\in{\mathbb{N}}\}$, which represents all positions in Single-delete Nim with three piles.  Let $N$ denote the set of all other positions.
 In the progression of this game, the total number of stones across all piles decreases monotonically. Since $P$ includes the terminal position $\la 1,1,1 \ra$, it is sufficient to show that no position in $P$ can be moved to any another position in $P$, and that any position in $N$ can be moved to some position in $P$.

First, assume $\la x,y,z \ra\in P$. In this case, since $v_{2}(x)=v_{2}(y)=v_{2}(z)$, it is sufficient to consider the next move where the pile with $z$ stones is eliminated and the pile with $x$ stones is split into two piles with $x'$ and $x''$ stones, resulting in the position $\la x',x'', y \ra$. By Proposition 2.1(i), if $v_{2}(x')=v_{2}(x'')$ then $v_{2}(x')<v_{2}(x)=v_{2}(y)$. This implies that $\la x',x'', y \ra$ cannot satisfy ($\bigstar$). Therefore, any position in $P$ cannot be moved to another position in $P$.

Next, assume $\la x,y,z \ra\in N$. Since $v_{2}(x)=v_{2}(y)=v_{2}(z)$ is not satisfied, it is sufficient to consider the case where $v_{2}(x)<v_{2}(y)$. By Proposition 2.2, the pile with $y$ stones can be split into two piles with $y'$ and $y''$ stones such that
$v_{2}(y')=v_{2}(y'')=v_{2}(x)$. Then,  by eliminating the pile with $z$ stones and splitting, the position can be moved to $\la x,y',y'' \ra\in P$. Thus, any position in $N$ can be moved to a position in $P$. $\Box$

\section{Single-delete Nim with four piles}
\subsection{Necessity of given conditions}
In this section, we present a detailed proof of Theorem 1.5, which establishes the necessary and sufficient conditions for a position in Single-delete Nim with four piles to qualify as a P-position.
We first briefly explain the necessity of some conditions (1) to (5) in Theorem 1.5 through a comparison of specific examples. The necessity of all the conditions will be discussed in detail in the proof of Lemma 3.2 later.

\begin{enumerate}
\renewcommand{\theenumi}{(\arabic{enumi})}
\item A specific example of a position in this case is $\la w,x,y,z \ra=\la 1440,864,$\\$672,1120 \ra$.
In its binary representation, it is
\[\begin{array}{ll}
w & 10110100000\\
x & 01101100000\\
y & 01010100000\\
z & 10001100000,
\end{array}
\]
where the rightmost 1's are in the same digits.
\item A specific example of a position in this case is $\la w,x,y,z \ra=\la 294,208,$\\$304,432 \ra$.
In its binary representation, it is
\[\begin{array}{ll}
w & 100100110\\
x & 011010000\\
y & 100110000\\
z & 110110000.
\end{array}
\]
In a similar position 
\[\begin{array}{ll}
w & 1001\underline{1}0110\\
x & 011010000\\
y & 100110000\\
z & 110110000
\end{array}
\]
where condition (2A) is not satisfied, 
by eliminating $x$ and partitioning $w$ into $w'$ and $w''$, 
it is possible to fulfill condition (2) in a move, resulting in 
\[\begin{array}{ll}
w' & 100100110\\
w''& 000010000\\
y & 100110000\\
z & 110110000.
\end{array}
\]
\item A specific example of a position in this case is $\la w,x,y,z \ra=\la 669,468,$\\$800,288 \ra$. 
In its binary representation, it is
\[\begin{array}{ll}
w & 1010011101\\
x & 0111010100\\
y & 1100100000\\
z & 0100100000.
\end{array}
\]
In a similar position 
\[\begin{array}{ll}
w & 10100\underline{0}1101\\
x & 01110\underline{0}0100\\
y & 1100100000\\
z & 0100100000
\end{array}
\]
where condition (3B) is not satisfied, 
by eliminating $z$ and partitioning $y$ into $y'$ and $y''$, 
it is possible to fulfill condition (3) in a move, resulting in 
\[\begin{array}{ll}
w & 1010001101\\
x & 0111000100\\
y' & 1100010000\\
y'' & 0000010000.
\end{array}
\]
\item  A specific example of a position in this case is $\la w,x,y,z \ra=\la 11133,$\\$12716,7136,13312 \ra$. 
In its binary representation, it is
\[\begin{array}{ll}
w & 10101101111101\\
x & 11000110101100\\
y & 01101111100000\\
z & 11010000000000.
\end{array}
\]
In a similar position 
\[\begin{array}{ll}
w & 10101101111101\\
x & 11000110101100\\
y & 011011\underline{0}1100000\\
z & 11010000000000
\end{array}
\]
where condition (4B) is not satisfied, 
by eliminating $x$ and partitioning $z$ into $z'$ and $z''$, 
it is possible to fulfill condition (3) in a move, resulting in 
\[\begin{array}{ll}
w & 10101101111101\\
y & 01101101100000\\
z' & 11001110000000\\
z'' & 00000010000000.
\end{array}
\]
\item A specific example of a position in this case is $\la w,x,y,z \ra=\la 45053,$\\$62932,32576,64512 \ra$. 
In its binary representation, it is
\[\begin{array}{ll}
w & 1010111111111101\\
x & 1111010111010100\\
y & 0111111101000000\\
z & 1111110000000000.
\end{array}
\]
In a similar position 
\[\begin{array}{ll}
w & 1010111111111101\\
x & 1111010111010100\\
y & 011\underline{0}111101000000\\
z & 1111110000000000
\end{array}
\]
where condition (5A) is not satisfied, 
by eliminating $x$ and partitioning $z$ into $z'$ and $z''$, 
it is possible to fulfill condition (4) in a move, resulting in 
\[\begin{array}{ll}
w & 1010111111111101\\
y & 0110111101000000\\
z' & 1110110000000000\\
z'' & 0001000000000000.
\end{array}
\]
\end{enumerate}

\subsection{Proof of Theorem 1.5}

In this section, we will prove Theorem 1.5. Let $G_{4}=\{\, \la w,x,y,z \ra \, | \, w,x,y,z$\\$\in\mathbb{N}\}$ be the set of all positions in Single-delete Nim with four piles. We say that $\la w,x,y,z \ra\in G_{4}$ is a {\it standard form} if it satisfies $v_{2}(w)\leq v_{2}(x)\leq v_{2}(y)\leq v_{2}(z)$. Rearranging the elements of $\la w,x,y,z \ra\in G_{4}$ in this order is referred to as the {\it standardization} of $\la w,x,y,z \ra$. 
We define the subset $P\subset G_{4}$ such that a position $\la w,x,y,z \ra$ belongs to $P$ if  and only if its standardization satisfies any of the conditions (1) to (5) in Theorem 1.5.  Let $N$ denote the set of all other positions in $G_4$. In other words, $P$ and $N$ correspond to the sets of P-positions and N-positions for this game, respectively.

In the progression of this game, the total sum of stones across all piles decreases monotonically, and since $P$ includes the terminal position $\la 1,1,1,1 \ra$, the proof of Theorem 1.5 requires demonstrating the following two lemmas.
\begin{lemma}
If $\la w,x,y,z \ra\in P$, then it  cannot be moved to any position in $P$.
\end{lemma}

\begin{lemma}
If $\la w,x,y,z \ra\in N$, then it can be moved to some position in $P$.
\end{lemma}
To prove these lemmas, we first define common notation. If the standard form of $\la w,x,y,z \ra$ satisfies condition ($i$) of Theorem 1.5 ($1\leq i \leq 5$), we denote $\la w,x,y,z \ra\in P_{i}$ and define the sets $P_{1},P_{2},P_{3},P_{4},P_{5}$ accordingly. Thus, $\displaystyle P=\bigcup_{i=1}^{5}P_{i}$, and by conditions (1) to (5) we have $P_{i}\cap P_{j}=\emptyset$  for any $i\neq j$.

We now present the following proposition to support our proof of the two lemmas.

\begin{proposition}
If  $w,x\in\mathbb{N}$ satisfy $I_{b+1}(w)=I_{b+1}(x)=1$ and $I_{c+1}(w)=I_{c+1}(x)$ ($b<c$), and additionally satisfy  $I_{j}(w)+I_{j}(x)\geq 1$ for all $j$ with $b+2\leq j \leq c$, then $I_{c+1}(w+x)=1$.
\end{proposition}

Considering the carry in binary addition, it is clear that the above proposition holds.
With this in mind, we will now proceed to prove Lemma 3.1 and Lemma 3.2.

\medskip

{\bf Proof of Lemma 3.1.} For a position $\la w,x,y,z \ra$ where $a=v_{2}(w), b=v_{2}(x), c=v_{2}(y), d=v_{2}(z)$,  
if  any of the following conditions hold: 
\begin{quote}
 ($\triangle$)\  $a=b=c<d$, $a=b<c=d$, $a=b<c<d$, $a<b=c<d$
\end{quote}
then none of the conditions (1) to (5) can be satisfied, and thus the position is not included in $P$.  
This fact will be frequently used in the explanation below. In this proof, when transitioning from $\la w,x,y,z \ra$ to $\la w',x',y',z' \ra$ as a move, we show that if  $\la w',x',y',z' \ra\in P$ then $\la w,x,y,z \ra\not\in P$. We consider the reverse of the operation:
\begin{itemize}
\item First, select two of $w',x',y',z'$ and combine their sum into a single number.
\item Then, add a new element $u'\in\mathbb{N}$.
\end{itemize}
We consider each case where $\la w',x',y',z' \ra\in P_{i}$ for $1\leq i \leq 5$. 
\begin{enumerate}
\renewcommand{\theenumi}{\alph{enumi})}
\item In case $\la w',x',y',z' \ra\in P_{1}$: \\
Since $v_{2}(w')=v_{2}(x')=v_{2}(y')=v_{2}(z')$, it is sufficient to consider only the case where the sum 
$y'+z'$ is formed in the reverse operation. In this situation, by Proposition 2.1(i), we have $v_{2}(w')=v_{2}(x')<v_{2}(y'+z')$. Therefore, regardless of the element $u'$ added, $\la u',w',x',y'+z' \ra$ corresponds to ($\triangle$), and cannot be included in any $P_{i}$.
\item In case $\la w',x',y',z' \ra\in P_{2}$: \\
Suppose $v_{2}(w')<v_{2}(x')=v_{2}(y')=v_{2}(z')$.
Let $b=v_{2}(x')$, and we will use the following facts: \\
($\clubsuit_{1}$) From (2A), we have $I_{b+1}(w')=0$.\\
($\clubsuit_{2}$) $I_{b+1}(w'+z')=1$?D
\begin{itemize}
\item We consider whether $\la u',w',x',y'+z' \ra$, formed by the reverse operation where the sum $y'+z'$ is taken and $u'$ is added, is included in $P$. Here, we have $v_{2}(w')<v_{2}(x')<v_{2}(y'+z')$. 
\begin{itemize}
\item If $v_{2}(u')>v_{2}(x')$, then by ($\clubsuit_{1}$), it cannot satisfy (3C), (4E) or (5F).
\item If $v_{2}(u')=v_{2}(x')$ or $v_{2}(u')=v_{2}(w')$, it corresponds to ($\triangle$).
\item If $v_{2}(u')<v_{2}(x')$ and $v_{2}(u')\neq v_{2}(w')$, then by ($\clubsuit_{1}$), it cannot satisfy (4C) or (5D).
\end{itemize}
Therefore, we can conclude that $\la u',w',x',y'+z' \ra\not\in P$.
\item We consider whether $\la u',w'+z',x',y' \ra$, formed by the reverse operation where the sum $w'+z'$ is taken and $u'$ is added, is included in $P$. Here, we have $v_{2}(w'+z')<v_{2}(x')=v_{2}(y')$.
\begin{itemize}
\item If $v_{2}(u')>v_{2}(y')$, it corresponds to ($\triangle$).
\item If $v_{2}(u')\leq v_{2}(y')$, then by ($\clubsuit_{2}$), it cannot satisfy (2A) or (3A).
\end{itemize}
Therefore, we can conclude that $\la u',w'+z',x',y' \ra\not\in P$.
\end{itemize}
\item In case $\la w',x',y',z' \ra\in P_{3}$: \\
Suppose $v_{2}(w')<v_{2}(x')<v_{2}(y')=v_{2}(z')$. 
Let $c=v_{2}(y')$, and we will use the following facts: \\
($\diamondsuit_{1}$) From (3A), we have $I_{c+1}(w')=I_{c+1}(x')=0$.\\ 
($\diamondsuit_{2}$) $I_{c+1}(x'+z')=1$.\\
($\diamondsuit_{3}$) From (3B) and Proposition 3.3, we have $I_{c+1}(w'+x')=1$.
\begin{itemize}
\item We consider whether $\la u',w',x',y'+z' \ra$, formed by the reverse operation where the sum $y'+z'$ is taken and $u'$ is added, is included in $P$. Here, we have $v_{2}(w')<v_{2}(x')<c<v_{2}(y'+z')$.
\begin{itemize}
\item If $v_{2}(u')>c$, then by ($\diamondsuit_{1}$), it cannot satisfy (3B), (4D) or (5E).
\item If $v_{2}(u')=c$, then by ($\diamondsuit_{1}$), it cannot satisfy (4C) or (5D).
\item If $v_{2}(u')<c$, $v_{2}(u')\neq v_{2}(w')$ and $v_{2}(u')\neq v_{2}(x')$, then by ($\diamondsuit_{1}$), it cannot satisfy (4B) or (5C).
\item If $v_{2}(u')=v_{2}(w')$ or $v_{2}(u')=v_{2}(x')$, it corresponds to ($\triangle$).
\end{itemize}
Therefore, we can conclude that $\la u',w',x',y'+z' \ra\not\in P$.
\item We consider whether $\la u',w',x'+z',y' \ra$, formed by the reverse operation where the sum $x'+z'$ is taken and $u'$ is added, is included in $P$. Here, we have $v_{2}(w')<v_{2}(x'+z')<v_{2}(y')$.
\begin{itemize}
\item If $v_{2}(u')>v_{2}(y')$, then by ($\diamondsuit_{1}$), it cannot satisfy (4C) or (5D).
\item If $v_{2}(u')=v_{2}(y')$, then by ($\diamondsuit_{2}$), it cannot satisfy (3A).
\item If $v_{2}(u')<v_{2}(y')$, $v_{2}(u')\neq v_{2}(x'+z')$ and $v_{2}(u')\neq v_{2}(w')$, then by ($\diamondsuit_{1}$) and ($\diamondsuit_{2}$), it cannot satisfy (4A) or (5B).
\item If $v_{2}(u')=v_{2}(x'+w')$ or $v_{2}(u')=v_{2}(z')$, it corresponds to ($\triangle$).
\end{itemize}
Therefore, we can conclude that $\la u',w',x'+z',y' \ra\not\in P$.
\item For $\la u',w'+z',x',y' \ra$, formed by the reverse operation where the sum $w'+z'$ is taken and $u'$ is added, 
 it can similarly be shown that it is not included in $P$, just as in the case where the sum $x'+z'$ is taken and $u'$ is added.
\item We consider whether $\la u',w'+x',y',z' \ra$, formed by the reverse operation where the sum $w'+x'$ is taken and $u'$ is added, is included in $P$. Here, we have $v_{2}(w'+x')<v_{2}(y')=v_{2}(z')$.
\begin{itemize}
\item If $v_{2}(u')>v_{2}(y')$, it corresponds to ($\triangle$).
\item If $v_{2}(u')\leq v_{2}(y')$ and $v_{2}(u')\neq v_{2}(w'+x')$, then by ($\diamondsuit_{3}$), it cannot satisfy (2A) or (3A).
\item If $v_{2}(u')=v_{2}(w'+x')$, it corresponds to ($\triangle$).
\end{itemize}
Therefore, we can conclude that $\la u',w'+x',y',z' \ra\not\in P$.
\end{itemize}
\item In case $\la w',x',y',z' \ra\in P_{4}$: \\
Suppose $v_{2}(w')<v_{2}(x')<v_{2}(y')<v_{2}(z')$. Let $c=v_{2}(y'), d=v_{2}(z')$,  and we will use the following facts: \\
($\heartsuit_{1}$) From (4A), we have $I_{d+1}(w')=I_{d+1}(x')=I_{d+1}(y')=0$.\\ 
($\heartsuit_{2}$) From (4C), we have $I_{c+1}(w')=I_{c+1}(x')=1$.\\ 
($\heartsuit_{3}$) From (4B), we have $I_{j}(w)+I_{j}(x)+I_{j}(y)\geq 2$ for $c+2\leq j \leq d$.\\
($\heartsuit_{4}$) From (4B) and by using similar arguments as Proposition 3.3, \\
\qquad \ we have $I_{d+1}(w'+x')=I_{d+1}(w'+y')=I_{d+1}(x'+y')=1$.
\begin{itemize}
\item We consider whether $\la u',w',x',y'+z' \ra$, formed by the reverse operation where the sum $y'+z'$ is taken and $u'$ is added, is included in $P$. Here, we have $v_{2}(w')<v_{2}(x')<v_{2}(y'+z')<d$.
\begin{itemize}
\item If $v_{2}(u')>d$, then by ($\heartsuit_{1}$), it cannot satisfy (4B) or (5C).
\item If $v_{2}(u')=d$, then by ($\heartsuit_{1}$), it cannot satisfy (4A) or (5B).
\item If $v_{2}(u')<d$, $v_{2}(u')\neq v_{2}(y'+z')$, $v_{2}(u')\neq v_{2}(x')$ and $v_{2}(u')\neq v_{2}(w')$, then by ($\heartsuit_{1}$) and ($\heartsuit_{3}$), it cannot satisfy (4A) or (5A).
\item If $v_{2}(u')=v_{2}(y'+z')$, then by ($\heartsuit_{2}$), it cannot satisfy (3A).
\item If $v_{2}(u')=v_{2}(x')$ or $v_{2}(u')=v_{2}(w')$, it corresponds to ($\triangle$).
\end{itemize}
Therefore, we can conclude that $\la u',w',x',y'+z' \ra\not\in P$.
\item For $\la u',w',x'+z',y' \ra$, formed by the reverse operation with the sum $x'+z'$ and addition of $u'$, and
 for $\la u',w'+z',x',y' \ra$, similarly with the sum $w'+z'$,  it can be shown that neither is included in $P$, just as in the case where the sum $y'+z'$ is taken and $u'$ is added.
\item We consider whether $\la u',w',x'+y',z' \ra$, formed by the reverse operation where the sum $x'+y'$ is taken and $u'$ is added, is included in $P$. Here, we have $v_{2}(w')<v_{2}(x'+y')<v_{2}(z')=d$.
\begin{itemize}
\item If $v_{2}(u')>v_{2}(z')$, then by ($\heartsuit_{1}$), it cannot satisfy (4C) or (5D).
\item If $v_{2}(u')=v_{2}(z')$, then by ($\heartsuit_{4}$), it cannot satisfy (3A).
\item If $v_{2}(u')<v_{2}(z')$, $v_{2}(u')\neq v_{2}(x'+y')$ and $v_{2}(u')\neq v_{2}(w')$, then by ($\heartsuit_{1}$) and ($\heartsuit_{4}$), it cannot satisfy (4A) or (5B).
\item If $v_{2}(u')=v_{2}(x'+y')$ or $v_{2}(u')=v_{2}(w')$, it corresponds to ($\triangle$).
\end{itemize}
Therefore, we can conclude that $\la u',w',x'+y',z' \ra\not\in P$.
\item For $\la u',w'+y',x',z' \ra$, formed by the reverse operation with the sum $w'+y'$ and addition of $u'$, and
 for  $\la u',w'+x',y',z' \ra$, similarly with the sum $w'+x'$,  it can be shown that neither is included in $P$, just as in the case where the sum $x'+y'$ is taken and $u'$ is added.
\end{itemize}
\item In case $\la w',x',y',z' \ra\in P_{5}$: \\
Suppose $v_{2}(w')<v_{2}(x')<v_{2}(y')<v_{2}(z')$. 
Let $b=v_{2}(x'), c=v_{2}(y'), d=v_{2}(z')$, and define the smallest $e$ such that $e>d$ and \\
($\spadesuit_{1}$) $I_{e+1}(w')=I_{e+1}(x')=I_{e+1}(y')=I_{e+1}(z')=0$.\\
We will use the following facts: \\
($\spadesuit_{2}$) From (5A), we have $I_{i}(w')+I_{i}(x')+I_{i}(y')+I_{i}(z')\geq 3$ for \\
\qquad $d+2\leq i\leq e$.\\
($\spadesuit_{3}$) From (5B), we have $I_{d+1}(w')=I_{d+1}(x')=I_{d+1}(y')=1$ and\\
\qquad \ $I_{d+1}(w'+z')=I_{d+1}(x'+z')=I_{d+1}(y'+z')=0$.\\
($\spadesuit_{4}$) From (5D), we have $I_{c+1}(w')=I_{c+1}(x')=1$.\\
($\spadesuit_{5}$) From (5C), we have $I_{j}(w)+I_{j}(x)+I_{j}(y)\geq 2$ for $c+2\leq j \leq d$, \\
\qquad and from (5E), 
$I_{k}(w)+I_{k}(x)\geq 1$ for $b+2\leq k \leq c$.\\
($\spadesuit_{6}$) From (5A) and by using similar arguments as Proposition 3.3, regarding the sum of any two of $w',x',y',z'$, $I_{e+1}(w'+x')=\cdots=I_{e+1}(y'+z')=1$.
\begin{itemize}
\item We consider whether $\la u',w',x',y'+z' \ra$, formed by the reverse operation where the sum $y'+z'$ is taken and $u'$ is added, is included in $P$. Here, we have $v_{2}(w')<v_{2}(x')<v_{2}(y'+z')=c$.
\begin{itemize}
\item If $v_{2}(u')>e$, then by ($\spadesuit_{1}$), it cannot satisfy (4B).
\item If $d<v_{2}(u')\leq e$, then by ($\spadesuit_{2}$) it cannot satisfy (4A) , and by ($\spadesuit_{1}$) and ($\spadesuit_{6}$) it cannot satisfy (5A).
\item If $v_{2}(u')=d$, then by ($\spadesuit_{3}$), it cannot satisfy (4A) and (5B).
\item If $v_{2}(u')<d$, $v_{2}(u')\neq v_{2}(y'+z')$, $v_{2}(u')\neq v_{2}(x')$ and  $v_{2}(u')\neq v_{2}(w')$, then by ($\spadesuit_{5}$) it cannot satisfy (4A), and by ($\spadesuit_{1}$) and ($\spadesuit_{6}$) it cannot satisfy (5A).
\item If $v_{2}(u')=v_{2}(y'+z')$, then by ($\spadesuit_{4}$), it cannot satisfy (3A).
\item If $v_{2}(u')=v_{2}(x')$ or $v_{2}(u')=v_{2}(w')$, it corresponds to ($\triangle$).
\end{itemize}
Therefore we can conclude that $\la u',w',x',y'+z' \ra\not\in P$.
\item For $\la u',w',x'+z',y' \ra$, formed by the reverse operation with the sum $x'+z'$ and addition of $u'$, and
 for $\la u',w'+z',x',y' \ra$, similarly with the sum $w'+z'$,  it can be shown that neither is included in $P$, just as in the case where the sum $y'+z'$ is taken and $u'$ is added.
\item We consider whether $\la u',w',x'+y',z' \ra$, formed by the reverse operation where the sum $x'+y'$ is taken and $u'$ is added, is included in $P$. Here, we have $v_{2}(w')<v_{2}(x'+y')<v_{2}(z')=d$.
\begin{itemize}
\item If $v_{2}(u')>e$, then by ($\spadesuit_{1}$), it cannot satisfy (4B).
\item If $v_{2}(z')<v_{2}(u')\leq e$, then by ($\spadesuit_{2}$) it cannot satisfy (4A), and by ($\spadesuit_{1}$) and ($\spadesuit_{6}$) it cannot satisfy (5A).
\item If $v_{2}(u')=v_{2}(z')$, then by ($\spadesuit_{3}$), it cannot satisfy (3A).
\item If $v_{2}(u')<v_{2}(z')$, $v_{2}(u')\neq v_{2}(x'+y')$ and $v_{2}(u')\neq v_{2}(w')$, then by ($\spadesuit_{5}$) it cannot satisfy (4A) and by ($\spadesuit_{1}$) and ($\spadesuit_{6}$) it cannot satisfy (5A).
\item If $v_{2}(u')=v_{2}(x'+y')$ or $v_{2}(u')=v_{2}(w')$, it corresponds to ($\triangle$).
\end{itemize}
Therefore, we can conclude that $\la u',w',x'+y',z' \ra\not\in P$.
\item For $\la u',w'+y',x',z' \ra$, formed by the reverse operation with the sum $w'+y'$ and addition of $u'$, and
 for $\la u',w'+x',y',z' \ra$, similarly with the sum $w'+x'$,  it can be shown that neither is included in $P$, just as in the case where the sum $x'+y'$ is taken and $u'$ is added. 
\end{itemize}
\end{enumerate}
Thus we completes the proof. $\Box$

\medskip

{\bf Proof of Lemma 3.2.} We explicitly show that for  $\la w,x,y,z \ra$ not included in $P$, it can be moved to some position that is included in $P$. Let $a=v_{2}(w), b=v_{2}(x), c=v_{2}(y), d=v_{2}(z)$.

\begin{enumerate}
\renewcommand{\theenumi}{\alph{enumi})}
\item In case any of the following holds: $a=b=c<d$, $a=b<c=d$, $a=b<c<d$ or, $a<b=c<d$:\\
First we assume $a=b<c<d$. By Proposition 2.2, there exists a partition $y=y'+y''$ such that $v_{2}(y')=v_{2}(y'')=a$ (For example, we set $y'=y-2^a$ and $y''=2^a$). Therefore, $\la w,x,y,z \ra$ can be moved to $\la w,x,y',y''\ra\in P_{1}$. Similarly, in the cases $a=b=c<d$, $ a=b<c=d$ and $a<b=c<d$, it can be moved to some position in $P_{1}$.
\item In case $a<b=c=d$ and (2A) is not satisfied: \\
By partitioning $w$ into $w'=w-2^{b}$ and $w''=2^{b}$, it can be moved to $\la w',w'',x,y\ra\in P_{2}$.
\item In case $a<b<c=d$ and at least one of (3A),(3B), or (3C) is not satisfied:
\begin{itemize}
\item Suppose that (3A) is not satisfied and $I_{c+1}(w)=1$.
By partitioning $w$ into $w'=w-2^{c}$ and $w''=2^{c}$, it can be moved to $\la w',w'',y,z\ra\in P_{2}$. The same applies in the case when $I_{c+1}(x)=1$.
\item Suppose that (3C) is not satisfied and $I_{b+1}(w)=0$. 
By partitioning $y$ into $y'=y-2^b$ and $y''=2^b$, it can be moved to $\la w,x,y',y''\ra\in P_{2}$.
\item Suppose that (3C) is satisfied but (3B) is not. 
For the smallest $k$ such that $c<k<d$ and $I_{k+1}(w)=I_{k+1}(x)=0$, by partitioning $y$ into $y'=y-2^k$ and $y''=2^k$, it can be moved to $\la w,x,y',y''\ra\in P_{3}$.
\end{itemize}
\item In case $a<b<c<d$ and (4A) is satisfied, and at least one of (4B), (4C), (4D), and (4E) is not satisfied:
\begin{itemize}
\item Suppose that (4E) is not satisfied, that is $I_{b+1}(w)=0$.
By partitioning $y$ into $y'=y-2^b$ and$y''=2^b$,  it can be moved to $\la w,x,y',y''\ra\in P_{2}$.
\item  Suppose that (4E) is satisfied but (4D) is not.
For the smallest $k$ such that $b<k<c$ and $I_{k+1}(w)=I_{k+1}(x)=0$, by partitioning $y$ into $y'=y-2^k$ and $y''=2^k$, it can be moved to $\la w,x,y',y''\ra\in P_{3}$.
\item Suppose that (4C) is not satisfied and $I_{c+1}(w)=0$.
By partitioning $z$ into $z'=z-2^{c}$ and $z''=2^{c}$, it can be moved to $\la w,y,z',z''\ra\in P_{2}$. The same applies in the case when $I_{c+1}(x)=0$.
\item Suppose (4C), (4D) and (4E) are satisfied but (4B) is not. 
Define the smallest $k$ such that $c<k<d$ and at least two of $I_{k+1}(w)$, $I_{k+1}(x)$, $I_{k+1}(y)$ are 0. For example, if $I_{k+1}(w)=I_{k+1}(x)=0$, then by partitioning $z$ into $z'=z-2^k$ and $z''=2^k$,  it can be moved to $\la w,x,z',z''\ra\in P_{3}$. The same applies in other cases.
\end{itemize}
\item In case $a<b<c<d$ and (4B), (4C), (4D) and (4E) are satisfied, but (4A) is not:

Since $I_{d+1}(w)+I_{d+1}(x)+I_{d+1}(y)\in\{1,2\}$, at least one of $I_{d+1}(w)$, $I_{d+1}(x)$ or $I_{d+1}(y)$ is 1, and at least one of them is 0. For example, if $I_{d+1}(w)=1$ and $I_{d+1}(x)=0$, then by partitioning $w$ into $w'=w-2^{d}$ and $w''=2^{d}$, it can be moved to $\la w',x,w'',z\ra\in P_{3}$. The same applies in other cases.
\item In case $a<b<c<d$ and (5B) is satisfied, and at least one of (5C), (5D), (5E) and (5F) is not satisfied:

It can be shown similarly as in case d) above.
\item In case $a<b<c<d$ and (5C), (5D), (5E), (5F) are satisfied, but (5B) is not:

It can be shown similarly as in case e) above.
\item In case $a<b<c<d$ and (5B), (5C), (5D), (5E), (5F) are satisfied, but (5A) is not:

Define the smallest $k$ such that $k>d$ and $I_{k+1}(w)+I_{k+1}(x)+I_{k+1}(y)+I_{k+1}(z)\in\{1,2\}$. Since 
at least one of $I_{k+1}(w)$, $I_{k+1}(x)$, $I_{k+1}(y)$ and $I_{k+1}(z)$ is 1, and at least two of them are 0.
 For example, if $I_{k+1}(w)=1$ and $ I_{k+1}(x)=I_{k+1}(y)=0$, then by partitioning $w$ into $w'=w-2^{k}$ and $w''=2^{k}$, it can be moved to $\la w',x,y,w''\ra\in P_{4}$. The same applies in other cases.
\end{enumerate}
Thus we completes the proof. $\Box$

\section{Single-delete Nim with five or more piles}

As discussed earlier, the winning and losing conditions for Single-delete Nim with four piles are significantly more complex than those for three or fewer piles. Additionally, it is not yet possible to even conjecture the conditions for scenarios with more piles. For Single-delete Nim with five or more piles, this paper is limited to presenting the following two almost self-evident propositions. In particular, Proposition 4.2 suggests that, unlike the conditions in Theorem 1.4 and Theorem 1.5(1), even if the number of times each pile's stone count can be divided by 2 is the same, the position does not necessarily become a P-position when there are more piles, indicating that the winning and losing conditions become even more complex.

\begin{proposition}
In Single-delete Nim with $n\geq 2$ piles, a position $\la x_{1},x_{2},$\\$x_{3},\ldots,x_{n} \ra$ is a P-position if all $x_{i}$ are odd. 
\end{proposition}
{\bf Proof. } Let $P_{0}$ be the set of positions where all $x_{i}$ are odd. 
When the player on their turn makes a move from $\la x_{1},x_{2},x_{3},\ldots,x_{n} \ra\in P_{0}$ to $\la x'_{1},x'_{2},x'_{3},\ldots,x'_{n} \ra$, exactly one of $x'_{1},x'_{2},x'_{3},\ldots,x'_{n}$ becomes even. 
In the next move, the opposite player eliminates one of the odd piles and splits the even pile into two odd piles, thus returning the position back to $P_{0}$. 
Therefore, at the position $\la x_{1},x_{2},x_{3},\ldots,x_{n} \ra\in P_{0}$, the player on their turn is always forced to take their  turn in a position belonging to $P_{0}$, eventually reaching the terminal position $\la 1,1,1,\ldots,1\ra$. $\Box$
\begin{proposition}
In Single-delete Nim with $n\geq 2$ piles, if the remainder when 
$n$ is divided by 3 is 2, then $\la 2,2,2,\ldots,2 \ra$ is a N-position. For any other value of $n$, $\la 2,2,2,\ldots,2 \ra$ is a P-position.
\end{proposition}

{\bf Proof. } In Single-delete Nim, starting from the initial position $\la 2,2,2,$\\$\ldots,2 \ra$, the positions that appear during the game can be described as having $k$ piles with 2 stones and $n-k$ piles with one stone ?i$1\leq k\leq n$?j. 
Considering that on each turn, the player can only remove a pile with 1 or 2 stones and split one of the piles with 2 stones, the game can be reduced to an equivalent stone-picking game of the following form by focusing on the value of $k$.
\begin{itemize}
\setlength{\leftskip}{0mm}
\item On the first turn, the player takes 2 stones from the total of $n$ stones. On subsequent turns, each player may take either 1 or 2 stones. 
\item The player who takes the last stone wins the game.
\end{itemize}
Therefore, when the remainder of $n$ divided by 3 is 2, the first player can always win by taking 2 stones initially. If the opponent then takes 2 stones, the first player should take 1 stone on the next turn; if the opponent takes 1 stone, the first player should take 2 stones. This strategy guarantees a win. The proof for the case when the remainder of $n$ divided by 3 is not 2 is omitted.
$\Box$

\section*{Acknowledgements}
We would like to express our gratitude to Dr. Koki Suetsugu, Prof. Tomoaki Abuku, and Prof. Ko Sakai for their invaluable comments regarding the study presented in this paper. 
This work was supported by JSPS KAKENHI Grant Number JP21K12191.

%
%
%
\bibliographystyle{splncs04}
\bibliography{SDNforarXiv}
%

%
%
%



\end{document}